\documentclass[a4paper]{article}
\usepackage{amsthm,amsmath,amssymb}
\newtheorem{teor}{Theorem}[section]
\newtheorem{defin}[teor]{Definition}
\newtheorem{lemm}[teor]{Lemma}
\newtheorem{osse}[teor]{Remark}
\newtheorem{prop}[teor]{Proposition}
\newtheorem{defi}[teor]{Definition}
\newtheorem{coro}[teor]{Corollary}
\newtheorem{prob}[teor]{Problem}
\newtheorem{assu}[teor]{Assumption}

\newcommand{\bele}{\begin{lemm}\begin{sl}}
\newcommand{\enle}{\end{sl}\end{lemm}}
\newcommand{\bedef}{\begin{defi}\begin{sl}}
\newcommand{\eddef}{\end{sl}\end{defi}}
\newcommand{\bete}{\begin{teor}\begin{sl}}
\newcommand{\ente}{\end{sl}\end{teor}}
\newcommand{\beos}{\begin{osse}\begin{rm}}
\newcommand{\eddos}{\end{rm}\end{osse}}
\newcommand{\bepr}{\begin{prop}\begin{sl}}
\newcommand{\empr}{\end{sl}\end{prop}}
\newcommand{\bepro}{\begin{prob}\begin{rm}}
\newcommand{\empro}{\end{rm}\end{prob}}
\newcommand{\bede}{\begin{defin}\begin{sl}}
\newcommand{\edde}{\end{sl}\end{defin}}
\newcommand{\beco}{\begin{coro}\begin{sl}}
\newcommand{\enco}{\end{sl}\end{coro}}

\newcommand{\beas}{\begin{assu}\begin{rm}}
\newcommand{\eddas}{\end{rm}\end{assu}}

\newcommand{\quext}{\quad\text}
\newcommand{\qquext}{\qquad\text}
\newcommand{\de}{\partial}

\newcommand{\RR}{\mathbb{R}}

\newcommand{\EE}{\mathbb{E}}

\newcommand{\NN}{\mathbb{N}}

\newcommand{\beeq}[1]{\begin{equation}\label{#1}}
\newcommand{\eddeq}{\end{equation}}

\newcommand{\beeqa}[1]{\begin{eqnarray}\label{#1}}
\newcommand{\eddeqa}{\end{eqnarray}}

\newcommand{\beal}[1]{\begin{align}\label{#1}}
\newcommand{\eddal}{\end{align}}

\newcommand{\bespl}[1]{\begin{split}\label{#1}}
\newcommand{\edspl}{\end{split}}

\newcommand{\bega}[1]{\begin{gather}\label{#1}}
\newcommand{\edga}{\end{gather}}

\newcommand{\beeqax}{\begin{eqnarray*}}
\newcommand{\eddeqax}{\end{eqnarray*}}

\def\qed{\ifmmode 
  \else \leavevmode\unskip\penalty9999 \hbox{}\nobreak\hfill
  \fi
  \quad\hbox{\hskip.5em\vrule width.4em height.6em depth.05em\hskip.1em}}
\def\endproofsym{\qed}
\renewenvironment{proof}[1][Proof]{\trivlist\item[\hskip\labelsep{\hskip0pt
    {\normalfont\scshape#1.}\hskip .321429\parindent}]\ignorespaces}
{\endproofsym\endtrivlist}
\def\endnobox{\def\endproofsym{}\end{proof}\def\endproofsym{\qed}}

\newcommand{\no}{\nonumber}

\newcommand{\beeqao}{\begin{eqnarray}\no}
\newcommand{\bealo}{\begin{align}\no}
\newcommand{\besplo}{\begin{split}\no}
\newcommand{\begao}{\begin{gather}\no}

\newcommand{\perogni}{\forall\,}
\newcommand{\esiste}{\exists\,}

\newcommand{\io}{\int_\Omega}

\newcommand{\epsi}{\varepsilon}

\newcommand{\dd}{_\delta}

\newcommand{\lhs}{left hand side}
\newcommand{\rhs}{right hand side}

\DeclareMathOperator{\ess}{ess}

\DeclareMathOperator{\deriv}{d}

\DeclareMathOperator{\dom}{dom}

\DeclareMathOperator{\sign}{sign}

\newcommand{\LDH}{L^2(0,T;H)}
\newcommand{\LDV}{L^2(0,T;V)}

\newcommand{\LIH}{L^\infty(0,T;H)}

\let\TeXchi\chi
\def\chi{{\setbox0 \hbox{\mathsurround0pt
$\TeXchi$}\hbox{\raise\dp0 \copy0 }}}

\newcommand{\teta}{\vartheta}

\newcommand{\calX}{{\mathcal X}}

\newcommand{\calE}{{\mathcal E}}

\newcommand{\calD}{{\mathcal D}}

\newcommand{\calJ}{{\mathcal J}}

\newcommand{\barfz}{\overline{f_0(\chi)}}

\newcommand{\Ov}{\overline}

\newcommand{\dit}{\deriv\!t}
\newcommand{\dis}{\deriv\!s}
\newcommand{\ditau}{\deriv\!\tau}

\newcommand{\diy}{\deriv\!y}
\newcommand{\dir}{\deriv\!r}

\newcommand{\ddt}{\frac{\deriv\!{}}{\dit}}

\newcommand{\calL}{{\cal L}}

\numberwithin{equation}{section}
\begin{document}

\title{Nonlocal phase-field
systems with general potentials}

\author{Maurizio Grasselli\\
Dipartimento di Matematica, Politecnico di Milano\\
Via Bonardi, 9\\
I-20133 Milano, Italy\\
E-mail: {\tt maurizio.grasselli@polimi.it}\\
\and
Giulio Schimperna\\
Dipartimento di Matematica, Universit\`a di Pavia,\\
Via Ferrata~1, I-27100 Pavia, Italy\\
E-mail: {\tt giusch04@unipv.it}
}

\maketitle
\begin{abstract}
\noindent We consider a phase-field model where the internal energy
depends on the order parameter $\chi$ in a nonlocal way. Therefore,
the resulting system consists of the energy balance equation coupled with a
nonlinear and nonlocal ODE for $\chi$. Such system has been analyzed
by several authors, in particular when the configuration potential is a smooth
double-well function. More recently, in the case of a potential defined on
$(-1,1)$ and singular at the endpoints,
the existence of a finite-dimensional global attractor has
been proven. Here we examine both the case of smooth potentials
as well as the case of physically realistic (e.g., logarithmic) singular
potentials. We prove well-posedness results and the eventual global boundedness
of solutions uniformly with respect to the initial data. In addition, we show
that the separation property holds in the case of singular potentials.
Thanks to these results, we are
able to demonstrate the existence of a finite-dimensional attractors in
the present cases as well.
\end{abstract}


\vspace{2mm}

\noindent {\bf AMS (MOS) subject clas\-si\-fi\-ca\-tion:}~~ 35B41,
35Q79, 37L30, 80A22.


%
%


\section{Introduction}
\label{sec:intro}

A well-known approach to study two-phase Stefan-like problems in more than one spatial dimension
is the so-called phase-field (or diffuse interface) method. Roughly speaking,
it consists in introducing an order parameter $\chi$ whose zero level set
substitutes for the sharp interface, while $\chi=\pm 1$ in the
higher/lower energy phases. The classical problem is thus
replaced by an order parameter dynamics, originated from the
study of critical phenomena, coupled with the energy balance
equation governing the temperature field. An important issue is to recover the
original interface conditions and this is usually done by (formal)
 asymptotic expansions. Moreover, diffuse interface models are quite
 effective from the numerical viewpoint since there
is no need for interface tracking. A significant and basic example
of phase-field system is due to G.~Caginalp (see \cite{C1}, cf.~also \cite{BS}), namely,
\begin{align}\label{calore0}
  & \teta_t + \frac{\ell}{2}\chi_t -k\Delta\teta = 0,\\
 \label{phase0}
  & \tau\chi_t -\xi\Delta\chi + W^\prime(\chi) = 2\teta,
\end{align}
on a given bounded domain $\Omega\subset \mathbb{R}^d$,
$d\in\{1,2,3\}$, for some time interval $(0,T)$. Here $\teta$ is a
rescaled temperature so that $\teta=0$ is the equilibrium melting
temperature, while $\ell$, $k$, $\tau$ and $\xi$ are given positive
constants  which represent the latent heat of fusion, the diffusivity, a
relaxation time and a correlation length, respectively. The function
$W$ is (the density of) potential energy associated with the phase
configuration. Such potential is usually a smooth double well function
(typically $W(r)=\frac{1}{8}(r^2-1)^2$). However, this is just a
convenient approximation of the physically more relevant
{\sl logarithmic potential}\/ generally taking the form
\begin{equation}\label{logpot}
  W(r) = (1+r)\log (1+r) + (1-r)\log(1-r) - \gamma r^2, \quad
   \gamma \ge 0.
\end{equation}
The mathematical literature on \eqref{calore0}-\eqref{phase0}
is rather vast and we confine
ourselves to quote the pioneering paper \cite{EZ} and the more
recent ones \cite{GMS,GPS1} (see also references therein).

In order to analyze the microscopic influences of anisotropy on the interface,
in \cite{C2} a phase-field system has been derived  from microscopic
considerations based on Statistical Mechanics. This system is similar
to \eqref{calore0}-\eqref{phase0} but for the diffusion term $\xi\Delta\chi$
which is replaced by $A:D^2\chi$, where $A\in \mathbb{R}^{d\times d}$ is
positive definite and $D^2\chi$ is the Hessian of $\chi$. However, this derivation
is performed by truncating the expansion of the interaction function
(see \cite[Prop.~2.4]{C2}, cf.~also \cite{CE} for higher-order approximations).
Then the author, by using formal asymptotics, deduces a modified Gibbs-Thompson
relation in 2D. More recently, by using the same procedure, a related phase-field
model has been obtained without approximating the interaction
function (see \cite{CCE1,CCE2}). Actually, working in a bounded domain
and choosing $\lambda=0$ in \cite[Sec.~2]{CCE2},
the system obtained there takes the following form:
\begin{align}\label{calore1}
  & \teta_t + \frac{\ell}{2}\chi_t -k\Delta\teta = 0,\\
 \label{phase1}
  & \varsigma^2\tau\chi_t - K_\varsigma * \chi + \kappa(x)\chi
  + W^\prime(\chi) = 2\varsigma\teta.
\end{align}
Here $K_\varsigma(x)=\varsigma^{-d}K(\varsigma^{-1} x)$ where
$\varsigma>0$ is an atomic length scale and $K: \RR^d \to \RR$ is a
sufficiently smooth interaction kernel satisfying $K(x)=K(-x)$ and
such that $\kappa(x):=\io K(x-y)\,\diy$ is bounded and nonnegative.

The asymptotic limit $\varsigma \searrow 0$ has been analyzed in
\cite{CCE2} and a new anisotropic interface condition has been obtained.
On the other hand, this class of systems was already
considered in some previous papers (cf., e.g.,
\cite{ABH,BC1,BCW,BHZ,FIRP} and their references). In particular,
rigorous mathematical results were proven for smooth potentials.
Well-posedness for $\Omega=\mathbb{R}$ and $d=1$ was
established in \cite{BCW} through semigroup theory. These results
were then extended to bounded $d$-dimensional domains with either
homogeneous Neumann or Dirichlet boundary conditions for $\teta$
(see \cite{BHZ,FIRP}). Regarding the longtime behavior, the
convergence of a solution to a single stationary state was shown in
\cite{FIRP} by means of a suitable nonsmooth version of the
{\L}ojasiewicz-Simon inequality. Existence of an absorbing set was
proven in \cite{BHZ} as well as an analysis of the $\omega$-limit
sets. More recently, the results of \cite{FIRP} have been extended to
to a class of singular unbounded potentials which does not include the
logarithmic ones (cf.~\cite{GPS2}). Actually, in \cite{GPS2} equation
\eqref{phase1} was modified by adding an inertial term of the form
$\alpha \chi_{tt}$, $\alpha>0$. However,
the results proved there also hold for $\alpha=0$. Then,
by exploiting \cite{GPS2}, the existence of a finite-dimensional global
attractor has been established in \cite{G}. Here we want to generalize
such results to both smooth potentials and more general singular
potentials (e.g., of the logarithmic type~\eqref{logpot}).
This goal is connected with the
property of the solutions of getting bounded in finite time uniformly
with respect to sufficiently general initial data. In addition, in the case
of singular potentials, a (uniform) separation property is also needed.
Here we prove all these properties for weak solutions originating from
initial data in the energy space. In particular, in the case of smooth
potentials, our results generalize the corresponding ones in
\cite{FIRP}. Moreover, for singular potentials, the separation property
holds instantaneously, namely for $t>0$, even though the initial datum
is a pure state (see Remark~\ref{waiting} below). As a consequence, we also demonstrate the existence of
a finite dimensional global attractor using the approach devised in
\cite{G}. This approach exploits the only source of compactness for
$\chi$, namely $K_\varsigma * \chi$. Note that we cannot expect smoothing
effects on $\chi$.

For the sake of simplicity, we choose the constants
in such a way that system
\eqref{calore1}-\eqref{phase1} can be rewritten in the form
\begin{align}\label{calore2}
  & \teta_t +  \chi_t -\Delta\teta = 0,\\
 \label{phase2}
  &\chi_t + J[\chi] + \kappa(x)\chi + W^\prime(\chi) = \teta,
\end{align}
in $\Omega\times(0,T)$. Here $J$ is a linear operator which is
a suitable generalization of the nonlocal convolution operator
introduced above (see Section~\ref{sec:main} below). Following~\cite{FIRP},
we endow the system with the following boundary and initial conditions
\begin{align}
\label{homDir}
  &\teta =0 \qquad\text{ on } \partial\Omega\times (0,T),\\
 \label{init0}
  & \teta|_{t=0}=\teta_0, \qquad \chi|_{t=0}=\chi_0 \qquad\text{ in }\Omega.
\end{align}
It is worth mentioning that there are also (mainly) existence and uniqueness
results for more refined nonlocal phase-field systems, formulated
with respect to the absolute temperature, which are
thermodynamically consistent also far from the equilibrium temperature
(cf.~\cite{KS1,KS2,KRS1,KRS2,SZ}).
It would be interesting and challenging to extend some of the present
results to such systems.

The plan of this paper goes as follows. The main results about
well-posedness and regularization properties of
the solutions are stated in Section~\ref{sec:main}.
Then, the corresponding theorems for smooth potentials
and singular potentials are proven in Section~\ref{sec:smooth}
and Section~\ref{sec:sing}, respectively. The final
Section~\ref{sec:attr} is devoted to the existence
of the global attractor.


\section{Well posedness and regularization results}
\label{sec:main}

In the sequel we will assume that $|\Omega|=1$, for simplicity. We
set $V:=H^1_0(\Omega)$, $H:=L^2(\Omega)$ and note by $\|
\cdot \|$ the norm in $H$, by $(\cdot,\cdot)$ the scalar product of
$H$, and by $\| \cdot \|_p$ the norm in $L^p(\Omega)$ for
$p\in[1,\infty]$. We also let $A$ be the Laplace operator with
homogeneous Dirichlet b.c., seen either as an unbounded linear
operator on $H$ with domain $V\cap H^2(\Omega)$ or as a bounded
linear operator from $V$ to its topological dual $V'$. Then
\eqref{calore2}-\eqref{homDir} can be rewritten as follows:
\begin{align}\label{calore}
  & \teta_t + \chi_t + A \teta = 0, \quext{in }\,V',\\
 \label{phase}
  & \chi_t + J[\chi] + f(\cdot,\chi) = \teta, \quext{in }\,H,
\end{align}
with
\begin{equation}\label{defif}
  f(x,r) = f_0(r) - \lambda(x) r,
  \end{equation}
where $\lambda\in L^\infty(\Omega)$ is a given function, and
\begin{equation}\label{defif0}
  f_0\in C^1(I;\RR), \qquad
   f_0'(r)\ge 0~~\perogni r \in I, \qquad
   f_0(0) = 0,
\end{equation}
where $I$, the {\sl domain}\/ of $f_0$, is an open and, possibly,
bounded interval of $\RR$ containing $0$. We also set
\begin{equation}\label{defiF}
  F_0(r):=\int_0^r f_0(s)\,\dis, \qquad
   F(x,r):=\int_0^r f(x,s)\,\dis = F_0(r) - \lambda(x) \frac{r^2}2.
\end{equation}
We assume $J$ be a linear operator such that
\begin{equation}\label{J1}
  J\in \calL(L^p(\Omega),L^p(\Omega)), \qquad
   \| J(u) \|_{L^p(\Omega)} \le L \| u \|_{L^p(\Omega)},
\end{equation}
for some $L>0$ independent of $p$ and all $p\in[1,+\infty]$.
Moreover we assume that
\begin{equation}\label{J2}
  \esiste p_*\in[1,+\infty):~~
  J\in \calL(L^{p_*}(\Omega),L^\infty(\Omega))
\end{equation}
and, finally,
\begin{equation}\label{Jcomp}
  J \quext{is a compact self-adjoint operator from $H$ to $H$}.
\end{equation}
Observe that the concrete form of the nonlocal operator $J$ (see \eqref{phase1})
satisfies assumptions \eqref{J1}-\eqref{Jcomp}, provided that $K$ is smooth enough.

We can then introduce the {\sl energy functional}
\begin{equation}\label{defiE}
  \calE(\teta,\chi):= \io \Big(
   \frac12 | \teta |^2 + F(\cdot,\chi)
     + \frac12 J[\chi]\chi \Big).
\end{equation}
It is immediate to realize that under the above assumptions
$\calE$ could be unbounded from below. Thus, we need some condition
implying that $\calE$ has some coercivity. In particular,
we will consider two different situations. The first one
deals with what we will call a {\sl smooth}\/ potential:
\beas\label{smooth:pot}
 We assume \eqref{defif0} with $I=\RR$. Moreover, we ask that
 \begin{equation}\label{Fsmooth}
   \kappa_f |r|^{1+\epsilon} - c_f
     \le f_0(r)\sign r
     \le C_f ( |r|^{1+\epsilon} +1 )
    \quad\perogni r\in \RR
 \end{equation}
 and for some $\epsilon>0$, $\kappa_f>0$, $c_f\ge 0$,
 $C_f>0$.
\eddas
%
%
\noindent%
We will speak, instead, of {\sl singular}\/ potentials
in the following case:
\beas\label{sing:pot}
 We assume \eqref{defif0} with $I$ an open and bounded interval
 of $\RR$ containing $0$. Moreover, we ask that
 \begin{equation}\label{Fsing}
   \lim_{r\to \de I} f_0(r)\sign r = +\infty.
 \end{equation}
\eddas
Our first result deals with the case when $F$ is smooth. In
this situation, we define the {\sl energy space}\/ (i.e., the
set of all $(\teta,\chi)$'s such that $\calE(\teta,\chi)$
is finite), as the Banach space
\begin{equation}\label{defiXsmooth}
  \calX := H\times L^{2+\epsilon}(\Omega).
\end{equation}
Indeed, by Assumption~\ref{smooth:pot}, $\calE(\teta,\chi)$ is finite
if and only if $(\teta,\chi)\in \calX$.
\bete\label{teo:smooth}
 Let\/ \eqref{J1}-\eqref{Jcomp} and\/
 {\rm Assumption~\ref{smooth:pot}}
 hold and let $(\teta_0,\chi_0)\in \calX$. Then, there exists
 one and only one couple $(\teta,\chi)$ satisfying,
 for all $T>0$,
 \begin{align}\label{regoteta}
   & \teta\in H^1(0,T;V') \cap \LIH \cap \LDV,\\
  \label{regochi}
   & \chi \in H^1(0,T;H) \cap L^\infty(0,T;L^{2+\epsilon}(\Omega)),
 \end{align}
 solving, a.e.~in~(0,T), system\/
 \eqref{calore}-\eqref{phase}, and enjoying
 the initial conditions
 \begin{equation}\label{init}
   \teta|_{t=0}=\teta_0, \quad
    \chi|_{t=0}=\chi_0,
     \quext{a.e.~in }\,\Omega.
 \end{equation}
 Moreover, there exist a time $T_0\ge 0$ depending
 on the ``initial energy'' $\EE_0:=\calE(\teta_0,\chi_0)$ and
 a constant $C_0$ independent of $\EE_0$,
 such that
 \begin{equation}\label{absorb}
   \| \teta(t) \|_{L^\infty(\Omega)}
    + \| \chi(t) \|_{L^\infty(\Omega)}
    + \| \chi_t(t) \|_{L^\infty(\Omega)}
    \le C_0 \quad\perogni t\ge T_0.
 \end{equation}
\ente
\noindent%
In the case when $F$ is singular, the
{\sl energy space}\/ is given by
\begin{equation}\label{defiXsing}
  \calX := \big \{ (\teta,\chi) \in H \times H:~
   F_0(\chi)\in L^1(\Omega) \big\}.
\end{equation}
Actually, since the domain $I$ of $F_0$ is bounded, it
is clear that $\calX\subset H\times L^\infty(\Omega)$.
Due to the constraint term $F_0$, $\calX$ is not a linear
space in this case. Nevertheless, it is
easy to prove that has a complete metric structure
with respect to a natural distance function
(see, e.g., \cite[Sec.~3]{RS} for details).
\bete\label{teo:sing}
 Let\/ \eqref{J1}-\eqref{Jcomp} and\/
 {\rm Assumption~\ref{sing:pot}}
 hold and let $(\teta_0,\chi_0)\in \calX$. Then, there exists
 one and only one couple $(\teta,\chi)$ satisfying,
 for all $T>0$,
 \begin{align}\label{regoteta2}
   & \teta\in H^1(0,T;V') \cap \LIH \cap \LDV,\\
  \label{regochi2}
   & \chi \in H^1(0,T;H), \quad F_0(\chi)\in L^\infty(0,T;L^1(\Omega)),\\
  \label{regofchi}
   & f_0(\chi)\in L^2(0,T;H),
 \end{align}
 and solving\/ \eqref{calore}-\eqref{phase} together with
 the initial conditions~\eqref{init}.
 Moreover, there exist a time $T_0\ge 0$ depending
 on the ``initial energy'' $\EE_0:=\calE(\teta_0,\chi_0)$ and
 a constant $C_0$ independent of $\EE_0$,
 such that\/ \eqref{absorb} holds, together with the\/
 {\rm separation property}
 \begin{equation}\label{separ}
   \| f(\chi) \|_{L^\infty(\Omega)}
    \le C_0 \quad\perogni t\ge T_0.
 \end{equation}
\ente
\noindent%
A couple $(\teta,\chi)$ in the condition either of
Theorem~\ref{teo:smooth} or of Theorem~\ref{teo:sing}
will be called an {\sl energy solution}\/ in what
follows.
\beos\label{source}
 An autonomous heat source term in equation \eqref{calore}
 could be easily handled (see, e.g., \cite{G}).
 Some care is required in the non-autonomous case
 (cf.~\cite{GPS1} for local systems).
\eddos


\section{Proof of Theorem~\ref{teo:smooth}}
\label{sec:smooth}

This existence proof is a slight generalization of \cite[Thm.~1.1]{FIRP}.
Thus we will proceed formally with the a priori estimates.
Uniqueness goes exactly as in \cite{FIRP},
while we will give all the details about \eqref{absorb}.
In the sequel we will note with the letter $c$
a generic positive constant, allowed to vary on
occurrence, depending only on $f_0$, $\lambda$ and $L$
(cf.~\eqref{J1}). In particular, we will always
assume $c$ to be independent of the initial data
and of time. The letter $\kappa$
will note the positive constants, depending on
the same quantities as $c$, appearing in estimates
from below.

\smallskip

\noindent%
{\bf Energy estimate.}~~%
We test \eqref{calore} by $\teta$, \eqref{phase} by
$\chi_t$ and take the sum. This gives (recall that $J$ is self-adjoint)
\begin{equation}\label{en-est}
  \ddt\calE(\teta,\chi)
   + \| \nabla \teta \|^2
   + \| \chi_t \|^2
  \le 0.
\end{equation}
Next, we test \eqref{phase} by $\chi$, to obtain
\begin{equation}\label{chi-est}
  \frac12 \ddt \| \chi \|^2
   + ( J[\chi],\chi )
   + ( f(\cdot,\chi),\chi )
   \le ( \teta,\chi ).
\end{equation}
By \eqref{J1}, Assumption~\ref{smooth:pot}, and $\lambda\in L^\infty(\Omega)$,
it is clear that
\begin{equation}\label{conto11}
  ( f(\cdot,\chi),\chi ) + ( J[\chi],\chi )
    \ge \kappa \io F(\cdot,\chi) + \kappa ( J[\chi],\chi ) - c
    \ge \kappa \| \chi \|_{2+\epsilon}^{2+\epsilon} + \kappa ( J[\chi],\chi ) - c.
\end{equation}
Moreover, by Poincar\'e's and Young's inequalities,
\begin{equation}\label{conto12}
   ( \teta,\chi )
    \le \frac12 \| \nabla \teta \|^2
     + c \| \chi \|^2
    \le \frac12 \| \nabla \teta \|^2
     + \sigma \| \chi \|_{2+\epsilon}^{2+\epsilon}
     + c_\sigma,
\end{equation}
for small $\sigma>0$ and $c_\sigma$ depending on $\sigma$.
Then, summing \eqref{en-est} with \eqref{chi-est} and taking
\eqref{conto11}, \eqref{conto12} into account,
we arrive at
\begin{equation}\label{diss-est}
  \ddt \Big( \calE(\teta,\chi)
   + \frac12 \| \chi \|^2 \Big)
   + \kappa \big( \calE(\teta,\chi)
   + \| \nabla \teta \|^2
   + \| \chi_t \|^2 \big)
  \le c.
\end{equation}
Integrating \eqref{diss-est}, we then obtain
\begin{equation}\label{diss-est-2}
  \calE(\teta(t),\chi(t)) \le c \big( \EE_0 e^{-\kappa t} + 1 \big),
\end{equation}
for some new value of $c$, independent of $T$.
In particular, by \eqref{Fsmooth}, this implies
\begin{equation}\label{diss-est-3}
  \| \teta(t) \|^2
   + \| \chi(t) \|_{2+\epsilon}^{2+\epsilon}
    \le c \big( \EE_0 e^{-\kappa t} + 1 \big).
\end{equation}
Moreover, integrating \eqref{en-est} over
the generic interval $(t,T)$,
and using \eqref{diss-est-2}, we infer
\begin{equation}\label{diss-integral}
  \int_t^T \big( \| \nabla \teta \|^2
   + \| \chi_t \|^2 \big)
   \le c \big( \EE_0 e^{-\kappa t} + 1 \big).
\end{equation}
Being $c$ independent of $T$,
the above bound can be rewritten
also for $T=+\infty$.

\smallskip

\noindent%
{\bf Regularization estimates for $\teta$.}~~%
From \eqref{diss-integral} and the Poincar\'e
inequality, it is clear that for any $s\in[0,+\infty)$
there exists $\tau=\tau(s)\in[s,s+1]$ such that
\begin{equation}\label{conto21}
  \| \teta(\tau) \|^2_V
   \le c \big( \EE_0 e^{-\kappa \tau} + 1 \big).
\end{equation}
Then, taking $s\ge 0$ and correspondingly choosing
$\tau\in[s,s+1]$ such that \eqref{conto21} holds,
testing \eqref{calore} by $-\Delta\teta$
and integrating over $(\tau,T)$ for a generic $T\ge \tau$,
recalling \eqref{diss-est}, and using H\"older's
and Young's inequalities, we obtain
\begin{equation}\label{conto22}
  \| \nabla\teta \|_{L^\infty(\tau,T;H)}^2
  + \int_\tau^T \| \Delta \teta \|^2
   \le c \big( \EE_0 e^{-\kappa \tau} + 1 \big).
\end{equation}
In particular, noting that the above holds at least
for some $\tau\in(0,1)$ (corresponding to the choice
$s=0$), we have
\begin{equation}\label{conto23}
  \| \teta \|_{L^\infty(t,T;V)}^2
   + \int_t^T \| \Delta \teta \|^2
   \le c \big( \EE_0 e^{-\kappa t} + 1 \big),
     \quad\perogni 1\le t\le T,
\end{equation}
still with $c$ independent of $T$.

\smallskip

\noindent%
{\bf Regularization estimates for $\chi$.}~~%
We will now work on a generic time interval
$(S,S+2)$ for $S\ge 1$. Then, as a consequence of estimate
\eqref{conto22} and interpolation, we have that
\begin{equation}\label{start:teta}
  \| \teta \|_{L^{10}((S,S+2)\times\Omega)}
   \le C.
\end{equation}
Here and in what follows, $C$ will always
denote a quantity of the form
\begin{equation}\label{defi:C}
  C = Q(\EE_0 e^{-\kappa S}),
\end{equation}
where $Q$ is a computable nonnegative-valued monotone
function, whose expression is allowed to vary
on occurrence, depending only on the fixed parameters of
the system.

That said, we choose a sequence of small time steps
$\tau_n$, $n\in \NN$, defined by
\begin{equation}\label{def:taun}
  \tau_n = \frac{3}{2\pi^2 n^2}
   \quext{so that }\,\sum_{n=1}^\infty \tau_n = \frac14
\end{equation}
and proceed by induction. Namely, we
set $t_0:=S$ and assume that, for $n\ge 1$,
given $t_{n-1}\ge S$, there
exists $t_n\in(t_{n-1},t_{n-1}+\tau_n)$ such that
\begin{equation}\label{conto31b}
  \| \chi(t_n) \|_{L^{2+n\epsilon}(\Omega)}^{2+n\epsilon}
   \le C \tau_n^{-1}.
\end{equation}
Notice that this is surely true as $n=1$ once one
sets $t_0=S$, thanks to \eqref{diss-est-3}.
Then, we test \eqref{phase} by
$|\chi|^{n\epsilon}\chi$.
In principle this would be not an admissible test
function (actually at this level \eqref{phase}
makes sense just as a relation in $L^2((0,T)\times \Omega)$
and $|\chi|^{n\epsilon}\chi$ needs not have
the $L^2$-summability).
However, the procedure could be easily justified by
using a truncation of $|\chi|^{n\epsilon}\chi$
as a test function, and then passing to the limit
w.r.t.~the truncation parameter via the monotone
convergence theorem (the details are left to
the reader).

Moreover, since we only need a finite number of iterations, we will not
take care of the dependence of the various constants
on $n$ and on $\epsilon$. We infer
\begin{equation}\label{conto32b}
  \frac1{2+n\epsilon} \ddt \| \chi \|^{2+n\epsilon}_{2+n\epsilon}
    + \big( f_0(\chi),|\chi|^{n\epsilon}\chi \big)
   \le \big( \teta, |\chi|^{n\epsilon}\chi \big)
    + \big( -J[\chi]+\lambda(\cdot)\chi, |\chi|^{n\epsilon}\chi \big).
\end{equation}
Then, by \eqref{Fsmooth},
\begin{equation}\label{conto33b}
   \big( f_0(\chi), |\chi|^{n\epsilon}\chi \big)
    \ge \kappa \| \chi \|_{2+(n+1)\epsilon}^{2+(n+1)\epsilon}
     - c.
\end{equation}
Moreover,
\begin{align}\no
  \big( \teta, |\chi|^{n \epsilon}\chi \big)
   & \le \big\| |\chi|^{n \epsilon}\chi \big\|_{\frac{2+(n+1)\epsilon}{1+n\epsilon}}
    \| \teta \|_{\frac{2+(n+1)\epsilon}{1+\epsilon}}\\
 \no
   & \le \| \chi \|_{2+(n+1)\epsilon}^{1+n\epsilon}
        \| \teta \|_{\frac{2+(n+1)\epsilon}{1+\epsilon}}\\
 \label{conto34b}
   & \le \frac\kappa4 \| \chi \|_{2+(n+1)\epsilon}^{2+(n+1)\epsilon}
    + c \| \teta \|_{\frac{2+(n+1)\epsilon}{1+\epsilon}}^{\frac{2+(n+1)\epsilon}{1+\epsilon}}.
\end{align}
Finally, by \eqref{J1} and $\lambda\in L^\infty(\Omega)$,
\begin{equation}\label{conto35b}
  \big(-J[\chi]+\lambda(\cdot)\chi,|\chi|^{n\epsilon}\chi \big)
   \le c \big\| |\chi|^{n \epsilon}\chi \big\|_{\frac{2+(n+1)\epsilon}{1+n\epsilon}}
    \| \chi \|_{\frac{2+(n+1)\epsilon}{1+\epsilon}}
   \le c \| \chi \|_{2+(n+1)\epsilon}^{2+n\epsilon}
   \le \frac\kappa4 \| \chi \|_{2+(n+1)\epsilon}^{2+(n+1)\epsilon}
    + c.
\end{equation}
Then, integrating \eqref{conto32b} over $(t_n,S+2)$,
recalling \eqref{start:teta},
and using \eqref{conto31b} and the induction hypothesis, we
arrive at
\begin{equation}\label{st32}
  \| \chi \|_{L^{\infty}(t_n,S+2;L^{2+n\epsilon}(\Omega))}^{2+n\epsilon}
   + \| \chi \|_{L^{2+(n+1)\epsilon}((t_n,S+2)\times \Omega)}^{2+(n+1)\epsilon}
   \le C + c \| \chi(t_n) \|_{L^{2+n\epsilon}(\Omega)}^{2+n\epsilon}
   \le C (1 + \tau_n^{-1}),
\end{equation}
where the second term on the \lhs\ ensures that
condition \eqref{conto31b}
will be fulfilled at the level $n+1$.
In particular, the procedure can be iterated
at least until $n$ satisfies the constraint
\begin{equation}\label{vincolon}
  \frac{2+(n+1)\epsilon}{1+\epsilon}\le 10,
   \quext{i.e., }\,n \le \frac{8+9\epsilon}{\epsilon}.
\end{equation}
More precisely, since $(8+9\epsilon)/\epsilon$ may not be an integer,
our method works at least until we reach some
$n_{\max}\ge 8(1+\epsilon^{-1})$.
Thus, we obtain (at least) the bound
\begin{equation}\label{st33}
  \| \chi \|_{L^{\infty}(S+1/4,S+2;L^{10+8\epsilon}(\Omega))}
   + \| \chi \|_{L^{10+9\epsilon}((S+1/4,S+2)\times \Omega)}
   \le C,
\end{equation}
where $C$ additionally depends on the chosen sequence $\tau_n$.
The {\sl upper}\/ bound in hypothesis~\eqref{Fsmooth} then gives also
\begin{equation}\label{st34}
  \| f(\cdot,\chi) \|_{L^9((S+1/4,S+2)\times \Omega)}
   \le C,
\end{equation}
whence, comparing terms in \eqref{phase} and using \eqref{start:teta},
assumption \eqref{J1} and \eqref{st33}-\eqref{st34}, we also obtain
\begin{equation}\label{st35}
  \| \chi_t \|_{L^9((S+1/4,S+2)\times \Omega)}
   \le C.
\end{equation}
With \eqref{st35} at disposal (indeed, any exponent $p>3$, in place of $9$,
would be sufficient for the purpose),
we can apply to equation \eqref{calore} (which is a linear PDE)
a standard Alikakos-Moser iteration method \cite{Al} with time-smoothing
(see, e.g., \cite[Lemma 3.5]{SSZ2}, cf.~also \cite[Chap.~III, Sec.~7]{LSU}),
to obtain
\begin{equation}\label{st36}
  \| \teta \|_{L^\infty((S+1/2,S+2)\times \Omega)}
   \le \Theta,
\end{equation}
for some $\Theta>0$ still having the form
$\Theta=Q(\EE_0 e^{-\kappa S})$.

\smallskip

Once we have \eqref{st36}, we can go back to \eqref{phase}
and perform further iterations of \eqref{conto32b},
restarting from $n=0$ with the choice of $t_0=S+1/2$.
Actually, \eqref{conto33b} can still
be used. On the other hand, thanks to \eqref{st36},
\eqref{conto34b} can be now modified this way:
\begin{equation}\label{conto34c}
  \big( \teta, |\chi|^{n \epsilon}\chi \big)
   \le \big\| |\chi|^{n \epsilon}\chi \big\|_{\frac{2+(n+1)\epsilon}{1+n\epsilon}}
    \| \teta \|_{\infty}
   \le \Theta \| \chi \|_{2+(n+1)\epsilon}^{1+n\epsilon}
   \le \frac\kappa4 \| \chi \|_{2+(n+1)\epsilon}^{2+(n+1)\epsilon}
    + c \Theta^{\frac{2+(n+1)\epsilon}{1+\epsilon}}.
\end{equation}
At this point, we can proceed with the iterations
exactly as before and notice
that, still in a {\sl finite}\/ number of steps, we arrive at
\begin{equation}\label{st41}
  \| \chi \|_{L^{\infty}(S+3/4,S+2;L^{p_*}(\Omega))}
   \le C.
\end{equation}
Consequently, using property \eqref{J2}, we also obtain
\begin{equation}\label{st42}
  \| J[\chi] \|_{L^{\infty}((S+3/4,S+2)\times \Omega))}
   \le M,
\end{equation}
for some constant $M>0$ having the same dependence of $C$.
\beos\label{no-moser}
 Unlike the case of parabolic type equations, it seems that,
 for equation \eqref{phase}, which has essentially an ODE structure,
 the Moser iteration method cannot be used to get directly an
 $L^\infty$-bound for $\chi$. Actually, if we try to
 iterate the above procedure {\sl infinitely many}\/ times,
 we readily notice that the constants appearing on the
 \rhs's of the estimates (cf., e.g., \eqref{conto34b}
 and \eqref{conto35b}) would explode.
 This is due to the fact that, while in true parabolic
 equations the Moser exponents
 grow {\sl exponentially}\/ with respect to $n$, in the present
 case, the growth is just {\sl linear}\/ (and, hence, too slow
 to take the constants under control). This fact also forces
 us to assume that the kernel $J$ has at least some
 regularizing effect (i.e., assumption \eqref{J2}).
\eddos
\noindent%
Once we have \eqref{st36} and \eqref{st42}
at our disposal, we can apply a
comparison principle to get the $L^\infty$-bound
for $\chi$. Namely, we have that, for $t\in(S+3/4,S+2)$
and a.e.~$x\in \Omega$, there holds
\begin{equation}\label{diff-eq}
  \chi_t(t,x) + f_0(\chi(t,x)) =
   \lambda(x)\chi(t,x) + \teta(t,x) - J[\chi](t,x).
\end{equation}
Let us now ``freeze'' the variable $x$. Actually,
thanks to the first \eqref{regochi2},
the map $t\mapsto \chi(t,x)$ is (Lipschitz) continuous
for a.e.~$x\in \Omega$.

Let now $\Lambda>0$ (to be chosen later) and
set
\begin{equation}\label{defiLambda+}
  \Lambda^+(x):=
   \big\{t\in(S+3/4,S+2):~\chi(t,x)\ge \Lambda\big\}.
\end{equation}
Then, using \eqref{Fsmooth}, \eqref{st36}, and \eqref{st42},
\eqref{diff-eq} gives
\begin{equation}\label{diff-ineq}
  \chi_t(t,x) + \kappa_f \chi(t,x)^{1+\epsilon}
   \le \lambda(x)\chi(t,x) + \Theta + M + c_f.
\end{equation}
for a.e.~$x\in\Omega$ and all $t \in \Lambda^+(x)$. Hence,
dividing by $\chi^{1+\epsilon}(t,x)$, we obtain
\begin{equation}\label{diff-ineq2}
  -\frac1\epsilon \ddt \chi^{-\epsilon}(t,x) + \kappa_f
   \le \lambda(x)\chi^{-\epsilon}(t,x)
     + (\Theta + M + c_f)\chi^{-(1+\epsilon)}(t,x).
\end{equation}
It is then clear that $\Lambda$ can be taken large enough
(in a way only depending on the $L^\infty$-norm of
$\lambda$ and on the known constants $\Theta$,
$M$, $c_f$ and $\kappa_f$) so that, for $t\in \Lambda^+(x)$,
\begin{equation}\label{diff-ineq3}
  -\frac1\epsilon \ddt \chi^{-\epsilon}(t,x) + \frac{\kappa_f}2
   \le 0.
\end{equation}
In particular, for such times $t$,
the function $t\mapsto \chi(t,x)$ is (strictly) decreasing.

This implies that, if $t\in (S+3/4,S+2)$
and $t\not\in \Lambda^+(x)$, then
$s\not\in \Lambda^+(x)$ for all $s\in [t,S+2)$.
In other words, if $\chi(t,x)$ is smaller than
$\Lambda$, it can never become larger than it.
Thus, assuming that $x$ is such that
$\chi_S(x):=\chi(S+3/4,x)\ge \Lambda$
(otherwise there is nothing to prove
in view of the preceding discussion),
integrating inequality \eqref{diff-ineq3}
over $(S+3/4,t)$, we obtain
\begin{equation}\label{diff-ineq4}
  \chi^{-\epsilon}(t,x)\ge
   \chi_S(x)^{-\epsilon} + \frac{\kappa_f\epsilon(t-S-3/4)}2,
   \quext{where }\,\chi_S(x):=\chi(S+3/4,x),
\end{equation}
at least for all $t\ge S$ such that $\chi(t,x) \ge \Lambda$.
Equivalently, we can write
\begin{equation}\label{diff-ineq5}
  \chi(t,x)\le
   \bigg( \frac{2\chi_S^{\epsilon}(x)}%
      {\kappa_f\epsilon\chi_S^\epsilon(x)(t-S-3/4)+2}
       \bigg)^{\frac1\epsilon}.
\end{equation}
Consequently, it is clear that there exists $\Lambda'>0$,
{\sl independent of the value of $\chi_S(x)$}, such that
\begin{equation}\label{diff-ineq6}
  \chi(t,x)\le \max\{\Lambda,\Lambda'\}
   \quext{for almost all }\, (t,x)\in (S+1,S+2)\times\Omega.
\end{equation}
For instance, one can take
$\Lambda' = (8\kappa_f^{-1}\epsilon^{-1})^{1/\epsilon}$.
Of course, a similar bound from below (of the form
$\chi(t,x)\ge -\max\{\Lambda,\Lambda'\}$) can be proved
in the same way. Thanks to the arbitrariness
of the starting time $S\in[1,+\infty)$
and recalling once more \eqref{st36}, we
have obtained the bounds for $\teta$ and $\chi$ in
\eqref{absorb}, for instance with the choice of $T_0=2$
(however, see Remark~\ref{waiting} below).
Then, the remaining bound for $\chi_t$ follows by
Assumption~\ref{smooth:pot}, \eqref{J1}, and
a comparison of terms in \eqref{phase}.
Theorem~\ref{teo:smooth} is proved.


\section{Proof of Theorem~\ref{teo:sing}}
\label{sec:sing}

This proof is presented in full detail (also for what concerns
existence) since, to the best of our knowledge, this
singular potential case has never been analyzed in the literature.


\subsection{Approximation}

We assume $I=\dom f_0=(-1,1)$ for simplicity.
We start by approximating the singular function $f_0$ by
a sequence $f\dd$, $\delta\in(0,1)$, such that
\begin{equation}\label{defifd}
  f\dd\in C^1(\RR;\RR), \qquad
   f\dd'(r)\ge 0~~\perogni r \in I, \qquad
   f\dd(0) = 0,
\end{equation}
for all $\delta\in(0,1)$. Moreover, we require that
\begin{equation}\label{Fdsmooth}
  \kappa_0 |r|^2 - c_0
    \le f\dd(r)\sign r
    \le C\dd ( |r|^2 +1 )
   \quad\perogni r\in \RR,~~\delta\in(0,1),
\end{equation}
where the constants $C\dd$ depend on $\delta$ (and, in fact,
will esplode as $\delta\searrow0$), while $\kappa_0$
and $c_0$ are assumed to be independent of $\delta$.
In other words, we are asking that
Assumption~\ref{smooth:pot} holds with $\epsilon=1$.
The singular character of $f_0$ ensures that the lower
bound in \eqref{Fdsmooth} can be assumed to
hold uniformly in~$\delta$.
Moreover, we assume that
\begin{align}\label{fmono}
  & f_{\delta_1}(r)\sign r \le f_{\delta_2}(r)\sign r
   \quad \perogni r\in I=(-1,1),~~
    \perogni 0 < \delta_2 < \delta_1 <1,\\
 \label{flim}
  & f\dd(r)\to f_0(r)\quext{uniformly on compact
   sets of }\,I=(-1,1),
\end{align}
whereas, for all $r\in\RR\setminus(-1,1)$,
$f\dd(r)\sign r\to +\infty$. We also set
\begin{equation}\label{defiFd}
  F\dd(r):=\int_0^r f\dd(s)\,\dis.
\end{equation}
Actually, it is easy to check that,
by \eqref{fmono}, $F_{\delta_1}\le F_{\delta_2}$ if
$\delta_2<\delta_1$. Moreover, we can assume that
\begin{equation}\label{propFd}
  f\dd(r)r - \lambda(x) r^2 \ge \kappa F\dd(r) - c
   \quad \perogni r\in\RR,~~\delta\in(0,1),
\end{equation}
for suitable constants $c\ge 0$, $\kappa>0$ independent of
$\delta$. The details of the construction of
this approximating family are left to the reader.
For instance, one possibility could be that of taking
$f\dd(r):= g\dd(r) + \delta^{-1} ((r-(1-\delta))^+)^2$
for $r\ge 0$ (and an analogous choice for $r<0$),
where $g\dd$ is the Yosida regularization of $f_0$
(cf., e.g., \cite{Br}).

We also notice that, if $\chi_0\in\calX$ (cf.~\eqref{defiXsing}),
recalling that $|\Omega|=1$, we have
\begin{equation}\label{unifiniz}
  \| \chi_0 \|_p \le 1 \quad\perogni
   p\in [1,\infty].
\end{equation}
Then, the above approximation permits to
apply Theorem~\ref{teo:smooth} to the system
\begin{align}\label{calored}
  & \teta_{\delta,t} + \chi_{\delta,t} + A \teta\dd = 0,\\
 \label{phased}
  & \chi_{\delta,t} + J[\chi\dd] + f\dd(\chi\dd) = \teta\dd + \lambda(x) \chi\dd,
\end{align}
with the initial conditions \eqref{init} (note that the initial
data are not approximated). This yields, for any $\delta\in(0,1)$,
a solution $(\teta\dd,\chi\dd)$ satisfying \eqref{regoteta}-\eqref{regochi}
(where $\epsilon=1$, cf.~\eqref{Fdsmooth}), together with
\eqref{absorb}. {\sl A priori}, these regularity properties
could depend on the approximation parameter $\delta$. However,
we shall see in a while that, in fact,
$\delta$-independent estimates are satisfied.


\subsection{Uniform estimates and passage to the limit}

In what follows we will assume that all constants
$\kappa$, $c$ are {\sl independent}\/ of $\delta$.
As we repeat the estimates performed in the proof
of Theorem~\ref{teo:smooth}, it is easy to realize
that, defining the approximate energy as
\begin{equation}\label{defiEd}
  \calE\dd(\teta,\chi):= \io \Big(
   \frac12 | \teta |^2
    + F\dd(\chi)
    - \frac{\lambda(\cdot)}2 \chi^2
     + \frac12 J[\chi]\chi \Big),
\end{equation}
the function $\calE\dd$ satisfies \eqref{diss-est}
with $\kappa$, $c$ independent of $\delta$.
Then, noting that
\begin{equation}\label{unifiniz2}
  \EE_{0,\delta}
   := \calE\dd(\teta_0,\chi_0)
   = \io \Big( \frac12 | \teta_0 |^2
       + F\dd(\chi_0) - \frac{\lambda(\cdot)}2 \chi_0^2
       + \frac12 J[\chi_0]\chi_0 \Big)
   \le \EE_0
\end{equation}
thanks to $F\dd\le F_0$, it is easy to check
that the ``Energy estimate'' of the previous
section can be repeated to obtain relations analogue to
\eqref{diss-est-2}, \eqref{diss-est-3},
\eqref{diss-integral}, that hold now uniformly
w.r.t.~$\delta$. Moreover, we can test \eqref{phased} by
$f\dd(\chi\dd)$ and use \eqref{diss-est-3}, \eqref{diss-integral},
and the properties of $J$ to infer
\begin{equation}\label{co41}
  \| f\dd(\chi\dd) \|_{L^2(t,t+1;H)} \le
   Q ( \EE_0 e^{-\kappa t} ),
   \quad \perogni t \ge 0,
\end{equation}
with $Q$ independent of $\delta$. We now show
that the estimates detailed
above suffice to take the limit $\delta\searrow 0$.
Actually, \eqref{diss-est-2}-\eqref{diss-integral}
and \eqref{co41} guarantee that, for any $T>0$,
\begin{align}\label{convteta}
  & \teta\dd\to \teta \quext{weakly star in }\,
   H^1(0,T;V') \cap \LIH \cap \LDV,\\
 \label{convchi}
  & \chi\dd\to \chi  \quext{weakly in }\,
   H^1(0,T;H),\\
 \label{convfchi}
  & f\dd(\chi\dd) \to \Ov{f_0(\chi)}  \quext{weakly in }\,
   \LDH.
\end{align}
Here and below, we adopt the convention of overlining
unidentified weak limits. Thanks to linearity and continuity
of operator $\calJ$, it is then easy to show that, at the
limit $\delta\searrow0$,
\begin{align}\label{calore-}
  & \teta_t + \chi_t + A \teta = 0,\\
 \label{phase-}
  & \chi_t + J[\chi] + \barfz - \lambda(x) \chi = \teta
\end{align}
and the initial conditions \eqref{init} are satisfied
as well.

Then, to conclude the proof, we have to show
the identification $\barfz=f_0(\chi)$
almost everywhere in $(0,T)\times\Omega$. To do this,
we follow with some variations the argument
given in \cite[Sec.~2.3]{FIRP}, which we report for the
reader's convenience.

First of all, letting
\begin{equation}\label{co11}
  \omega\dd:[0,T] \to \RR, \qquad
  \omega\dd(t):= \| \chi\dd(t) - \chi(t) \|^2,
\end{equation}
using \eqref{convchi} it is a standard check to verify that
\begin{equation}\label{co12}
  \| \omega\dd \|_{H^1(0,T)} \le c.
\end{equation}
Hence, we can assume that (here and below, all convergence
relations are intended up to extraction of non-relabelled
subsequences)
\begin{equation}\label{co13}
  \omega\dd \to \omega \quext{strongly in }\,C^0([0,T]),
\end{equation}
where $\omega$ is continuous and nonnegative.
Now, as a further consequence of \eqref{convchi},
we have that
\begin{equation}\label{st11b}
  \chi\dd \to \chi \quext{in }\, C_w([0,T];H).
\end{equation}
In particular, for {\sl all}\/ $t\in[0,T]$,
$\chi\dd(t)$ converges to $\chi(t)$ weakly in $H$.
%
%
%
%
Next, we compute the difference
between \eqref{phased} and \eqref{phase-},
test it $\chi\dd - \chi$, and integrate with respect to
the space variables. This gives
\begin{equation}\label{co15}
  \frac12 \ddt \| \chi\dd - \chi \|^2
   + \big( f\dd(\chi\dd) - \barfz, \chi\dd - \chi \big)
   = \big( \teta\dd - J[\chi\dd] + \lambda(\cdot)\chi\dd
    - \teta + J[\chi] - \lambda(\cdot)\chi, \chi\dd - \chi \big).
\end{equation}
Now, we notice that, by \eqref{co41}, \eqref{convchi} and
the first inequality in \eqref{Fdsmooth},
\begin{equation}\label{co16}
  \big\| f\dd(\chi\dd) \chi\dd \|_{L^{4/3}((0,T)\times \Omega)} \le c.
\end{equation}
Consequently,
\begin{equation}\label{co17}
  f\dd(\chi\dd) \chi\dd  \to \Ov{f_0(\chi) \chi}
   \quext{weakly in }\,L^{4/3}((0,T)\times \Omega).
\end{equation}
By definition of subdifferential, we have,
almost everywhere in $(0,T)\times \Omega$,
\begin{equation}\label{co18}
  f\dd(\chi\dd)(\chi\dd - \chi)
   \ge F\dd(\chi\dd) - F\dd(\chi)
   \ge F\dd(\chi\dd) - F_0(\chi).
\end{equation}
Let us now test \eqref{co18} by a nonnegative
test function $\phi \in \calD((0,T)\times \Omega)$
and integrate. Then,
by convexity and lower semicontinuity of the functional
\begin{equation}\label{co19}
  L^2((0,T)\times \Omega) \to \RR, \qquad
   v \mapsto \iint_{(0,T)\times \Omega}
    F_0(v) \phi,
\end{equation}
using \eqref{convfchi} and \eqref{co17}, we obtain that
\begin{align}\no
  & \iint_{(0,T)\times \Omega}
    \Ov{f_0(\chi) \chi} \phi
   - \iint_{(0,T)\times \Omega} \Ov{f_0(\chi)} \chi \phi\\
 \no
  & \mbox{}~~~~~
  = \lim_{\delta\searrow0}
   \iint_{(0,T)\times \Omega}
   f\dd(\chi\dd) (\chi\dd - \chi) \phi \\
 \label{co20}
  & \mbox{}~~~~~
  \ge \liminf_{\delta\searrow0}
   \iint_{(0,T)\times \Omega}
   \big( F\dd(\chi\dd) - F_0(\chi) \big) \phi
  \ge 0.
\end{align}
To deduce the last inequality we have used the
fact that the family of functionals $\{F\dd\}$,
being monotone increasing with respect to $\delta$ going
to 0, converges to $F_0$ in the sense of Mosco
(see, e.g., \cite{At}) in $L^2((0,T)\times \Omega)$.
In particular, we used here the $\liminf$-property
of Mosco-convergence:
\begin{equation}\label{Mosco}
  \iint_{(0,T)\times \Omega}  F_0(v)\phi
    \le \liminf_{\delta\searrow0} \iint_{(0,T)\times \Omega}  F\dd(v\dd) \phi
   \quext{for all }\, v\dd \to v~~\text{weakly in }\,
     L^2((0,T)\times \Omega).
\end{equation}
Thus, we have, almost everywhere in $(0,T)\times \Omega$,
\begin{equation}\label{co21}
  \Ov{f_0(\chi) \chi} \ge \Ov{f_0(\chi)}\chi.
\end{equation}
Moreover, we notice that, thanks to \eqref{convteta},
\eqref{convchi}, the Aubin-Lions Lemma, and
assumption \eqref{Jcomp},
\begin{equation}\label{co21b}
  \big( \barfz + \teta\dd - J[\chi\dd] - \teta + J[\chi], \chi\dd - \chi \big)
   \to 0,
\end{equation}
at least in the sense of distributions over $(0,T)$.

Then, we can take the limit, as $\delta\searrow 0$,
of \eqref{co15}. Using \eqref{co17} and
\eqref{co21}, and noting that the time-derivative operator
is linear and continuous with respect to distributional
convergence, we then obtain
\begin{equation}\label{co23}
  \frac12 \ddt \lim_{\delta\searrow 0} \| \chi\dd - \chi \|^2
   \le c_\lambda \lim_{\delta\searrow 0} \| \chi\dd - \chi \|^2,
\end{equation}
or, equivalently,
\begin{equation}\label{co24}
  \frac12 \ddt \omega(t)
   \le c_\lambda \omega(t).
\end{equation}
Since $\omega$ is nonnegative
and $\omega(0)=0$, we then obtain $\omega(t)=0$
for all $t\in[0,T]$. In other words,
\begin{equation}\label{co25}
  \chi\dd(t) \to \chi(t) \quext{strongly in }\,H,
   \quad\perogni t\in[0,T].
\end{equation}
This fact, combined with \eqref{convchi} gives
\begin{equation}\label{co27}
  \chi\dd \to \chi \quext{(at least) strongly in }\,
  {L^2(0,T;H)},
\end{equation}
which entails in particular
$\Ov{f_0(\chi)}=f_0(\chi)$. Thus, \eqref{phase-}
reduces to \eqref{phase}, as desired.


\subsection{Regularization estimates and separation property}

The above procedure is sufficient to get existence
of an energy solution to our system under
Assumption~\ref{sing:pot}. Uniqueness of this solution is
proved as in the other cases.

To prove rigorously the separation property \eqref{separ}, we
go back to the $\delta$-system \eqref{calored}-\eqref{phased}
and start by noticing that the regularization estimates of
Section~\ref{sec:smooth} hold {\sl uniformly}\/ in $\delta$.
Actually, concerning the ``Regularization estimates for $\teta$'' it is
easy to see that nothing changes and the analogue of
\eqref{conto21}-\eqref{conto23} hold uniformly in $\delta$.
Concerning the ``Regularization estimates for $\chi$'', we can
proceed as before (where we now have $\epsilon=1$, of course), until
we reach estimate \eqref{st33}. Indeed, in this part of the Moser
iteration, we only use the estimates
\eqref{diss-est-3}-\eqref{diss-integral} and
\eqref{conto21}-\eqref{conto23}, which are uniform in $\delta$, and
the estimate from below (i.e., the first inequality) in
\eqref{Fdsmooth}, which is also independent of $\delta$.
Consequently, we now have the analogue of \eqref{st33}, which, for
$\epsilon = 1$ and in the current notation, becomes
\begin{equation}\label{st33d}
  \| \chi\dd \|_{L^{\infty}(S+1/4,S+2;L^{18}(\Omega))}
   + \| \chi\dd \|_{L^{19}((S+1/4,S+2)\times \Omega)}
   \le C.
\end{equation}
Here and in what follows, all constants $c$, $\kappa$ and $C$ have the
same meaning as in the previous section and, in addition,
are assumed to be independent of $\delta$.

However, at this point we can no longer deduce the analogue of
\eqref{st34} directly, since this requires use
of the {\sl upper}\/ bound in \eqref{Fdsmooth},
where the constants {\sl do depend}\/ on $\delta$.

We then have to proceed with some more care and
set, for $p\in [1,\infty)$,
\begin{equation}\label{defiphip}
  \phi\dd^p(s) := \int_0^s |f\dd(r)|^{p} \sign r \,\dir
   \le |f\dd(s)|^{p} |s|.
\end{equation}
Then, we test \eqref{phased} by $|f\dd(\chi\dd)|^{p} \sign \chi\dd$,
with $p$ to be chosen below. This gives
\begin{align}\no
  \ddt \io \phi\dd^p(\chi\dd)
   + \| f\dd(\chi\dd) \|_{p+1}^{p+1}
  & = \big( |f\dd(\chi\dd)|^{p} \sign \chi\dd, \teta\dd + \lambda(\cdot) \chi\dd - J[\chi\dd] \big)\\
 \label{co31}
  & \le \frac12 \| f\dd(\chi\dd) \|_{p+1}^{p+1}
  + c \big( \| \teta\dd \|_{p+1}^{p+1} + \| \chi\dd \|_{p+1}^{p+1} \big).
\end{align}
Thus, integrating in time over $(\tau,t)$, where $S+1/4\le \tau \le t
\le S+2$, we arrive at
\begin{align}\no
  & \io \phi\dd^p(\chi\dd(t))
   + \frac12 \| f\dd(\chi\dd) \|^{p+1}_{L^{p+1}((\tau,t)\times \Omega)}\\
 \label{co32}
  & \mbox{}~~~~~
  \le \io \phi\dd^p(\chi\dd(\tau))
  + c \big( \| \teta\dd \|^{p+1}_{L^{p+1}((\tau,t)\times \Omega)}
   + \| \chi\dd \|^{p+1}_{L^{p+1}((\tau,t)\times \Omega)} \big).
\end{align}
Let us now first choose $p=5/3$. Then, according to \eqref{co41},
we can take $\tau\in (S+1/4,S+5/16)$ such that
\begin{equation}\label{co41b}
  \| f\dd(\chi\dd(\tau)) \|_2^2 \le C.
\end{equation}
Then, by the inequality in \eqref{defiphip},
\begin{equation}\label{co41c}
  \io \phi\dd^{5/3}(\chi\dd(\tau))
   \le \big\| |f\dd(\chi\dd(\tau))|^{5/3} \big\|_{6/5}
    \big\| \chi\dd(\tau) \big\|_6
   \le \big\| f\dd(\chi\dd(\tau)) \big\|_2^{5/3}
    \big\| \chi\dd(\tau) \big\|_6
\end{equation}
and both factors on the \rhs\ are controlled, uniformly in $\delta$,
thanks to \eqref{co41b} and \eqref{st33d}, respectively.
Then, noting that the other terms on the \rhs\ of \eqref{co32}
are estimated due to \eqref{start:teta} and \eqref{st33d},
we readily infer
\begin{equation}\label{co41d}
   \| f\dd(\chi\dd) \|_{L^{8/3}((S+5/16,S+2)\times\Omega)} \le C.
\end{equation}
At this point, we iterate the procedure by
choosing now $p=7/3$ in \eqref{co32}, taking $\tau\in (S+5/16,S+3/8)$
such that the analogue of \eqref{co41b} with
$8/3$ in place of $2$ holds (this is possible
thanks to \eqref{co41d}),
and replacing \eqref{co41c} with
\begin{equation}\label{co41c2}
  \io \phi\dd^{7/3}(\chi\dd(\tau))
   \le \big\| |f\dd(\chi\dd(\tau))|^{7/3} \big\|_{8/7}
    \big\| \chi\dd(\tau) \big\|_8
   \le \big\| f\dd(\chi\dd(\tau)) \big\|_{8/3}^{7/3}
    \big\| \chi\dd(\tau) \big\|_8.
\end{equation}
Thus, we finally arrive at
\begin{equation}\label{co41e}
   \| f\dd(\chi\dd) \|_{L^{10/3}((S+3/8,S+2)\times\Omega)} \le C.
\end{equation}
A comparison of terms in \eqref{phased}
permits now to get (in place of \eqref{st35})
\begin{equation}\label{co42}
   \| \chi_{\delta,t} \|_{L^{10/3}((S+3/8,S+2)\times\Omega)} \le C
\end{equation}
and the exponent $10/3$ is still sufficient to
operate the Moser iteration argument with smoothing
for $\teta$. Thus, we arrive also in this case
at the analogue of relation \eqref{st36}, with
$\Theta$ independent of $\delta$.

With this relation at disposal, we can repeat the
ODE argument of Section~\ref{sec:smooth}, with essentially
no variation. Actually, it is sufficient to use the lower
bound in \eqref{Fdsmooth}, which is uniform in $\delta$.
Summarizing, we have obtained
\begin{equation}\label{ThetaM}
  \| \teta\dd \|_{L^\infty((1,\infty)\times\Omega)} \le \Theta, \qquad
   \| \chi\dd \|_{L^\infty((1,\infty)\times\Omega)} \le M,
\end{equation}
with constants $\Theta$ and $M$ independent of $\delta$.
\beos\label{relevant}
 It is worth noting that, at the limit step, we have
 for free that $-1\le \chi \le 1$ almost everywhere
 (and starting from $t=0$), since $F_0$ is singular.
 However, \eqref{ThetaM}
 says something more, i.e., that we have,
 for $t\ge 1$, a uniform $L^\infty$-bound
 independent of the approximation parameter
 for {\sl both}\/ components of the approximate solution.
 This is a nontrivial information especially as far as
 $\teta\dd$ is concerned (indeed, we just know
 that the initial datum $\teta_0$ lies in $H$).
\eddos
\noindent%
As a consequence of \eqref{ThetaM},
there exists (a new) $\Lambda>0$ depending only
on $\lambda$, $J$, $\Theta$ and $M$, such that
\begin{equation}\label{compar1}
  \chi_{\delta,t}(t,x) + f\dd(\chi\dd(t,x))
   \le \Lambda \quext{for a.e.~}\, (t,x)\in (1,\infty)\times\Omega.
\end{equation}
Moreover,
\begin{equation}\label{compar2}
  \chi\dd(1,x) \le M
   \quext{a.e.~in }\, \Omega.
\end{equation}
Let now $\epsi\in(0,1)$ such that $f_0(1-\epsi)=3\Lambda$
(cf.~\eqref{Fsing}).
Then, there exists $\delta_0\in(0,1)$ such that
$f\dd(1-\epsi)\ge 2\Lambda$ for all $\delta\in[0,\delta_0]$.
In particular, by monotonicity (cf.~\eqref{fmono})
we have that $f\dd(r)\ge 2\Lambda$
for all $r\ge 1-\epsi$ and $\delta\in(0,\delta_0]$.
We now claim that there exists a time $T_0$ independent
of $\delta$ such that, for (almost) all $x\in\Omega$,
all $t\ge T_0$, and all $\delta\in(0,\delta_0]$, there
holds that $\chi\dd(t,x)\le 1-\epsi$.
Indeed, thanks to the above discussion, we have that,
if for some $\delta\in(0,\delta_0]$, $t\ge 1$,
and $x\in \Omega$, it is $\chi\dd(t,x) > 1-\epsi$,
then $t\mapsto\chi\dd(t,x)$ is decreasing. Thus, once
$t\mapsto \chi\dd(t,x)$ enters the region where
$\chi\dd(t,x)\le 1-\epsi$, it never exits from it.
Thus, freezing $x$ as before and
assuming that $\chi\dd(1,x) > 1-\epsi$
(otherwise there is nothing to prove), from
\eqref{compar1} we have that
\begin{equation}\label{compar4}
  \chi_{\delta,t}(t,x) \le -\Lambda,
\end{equation}
at least as long as $\chi(t,x)$
remains larger than $1-\epsi$.
Then, integrating \eqref{compar1} in time
and using \eqref{compar2}, we infer
\begin{equation}\label{compar5}
  \chi_{\delta}(t,x) \le \chi\dd(1,x)
   - \Lambda (t-1) \le M - \Lambda(t-1),
\end{equation}
whence
\begin{equation}\label{compar6}
  \chi_\delta(t,x) \le 1-\epsi
   \quext{for all }\,t\ge \frac{M-(1-\epsi)}\Lambda+1
\end{equation}
and this bound is uniform with respect to $\delta$. Proving
the lower bound in the same way and passing to the limit
$\delta\searrow0$, we finally obtain~\eqref{separ}, as desired.
\beos\label{waiting}
 Looking at the statement (and at the proof) of
 \eqref{separ} one could think that the
 separation property occurs only after some
 waiting time. However, the property was given in that
 form just for the sake of simplicity. Indeed,
 refining a bit the arguments in the proof, it
 is easy to demonstrate that \eqref{separ} is in fact
 an instantaneous property. Namely, there holds
 \begin{equation}\label{separepsi}
    -1+\epsi(\tau) \le \chi(t,x)
     \le 1-\epsi(\tau) \quext{for all }
     (t,x)\in (\tau,+\infty)\times\Omega
     ~~\text{and all }\,\tau>0,
 \end{equation}
 where $\epsi(\tau)$ goes to $0$ as $\tau\searrow0$.
 The details are left to the reader.
\eddos


\section{Global attractors}
\label{sec:attr}

Here we establish the existence of a finite-dimensional global
attractor for both smooth and singular potentials. The
technique is the same as the one used
in~\cite[Proof of Thm.~4.1]{G} for singular
unbounded potentials. However, here we start from
very general initial data in the energy space and we exploit
the previous regularization results to define a semigroup acting
on a convenient invariant set.

Let us consider the case of smooth potentials first. From
Theorem~\ref{teo:smooth} it is clear that we can define a semigroup
$S(t): \calX \to \calX$ (cf.~\eqref{defiXsmooth}) by setting
$(\teta(t),\chi(t)):= S(t)(\teta_0,\chi_0)$, where $(\teta,\chi)$ is
the unique (energy) solution to \eqref{calore}-\eqref{phase},
\eqref{init}. Note that $\teta \in C^0([0,+\infty);H)$ while $\chi\in
C^0([0,+\infty);H) \cap
C^0_w([0,+\infty);L^{2+\epsilon}(\Omega))$. Moreover, on
account of the Lipschitz continuity estimate \cite[(1.11)]{FIRP}, this
semigroup is also closed in the sense of \cite{PZ}. Thanks to
\eqref{absorb} and \eqref{conto23}, $S(t)$ has an absorbing set
$\mathcal{B}$ which is bounded in
$(V\cap L^\infty(\Omega))\times L^\infty(\Omega)$. Hence we can find
$t_0\geq 0$ such that $S(t)\mathcal{B} \subset
\mathcal{B}$ for all $t\geq t_0$. Then, without loss of generality, we
can suppose that $t_0=0$ and assume that $\mathcal{B}$ is an
invariant set for $S(t)$. Moreover, we can endow $\mathcal{B}$ with
the $V\times H$-metric and obtain a complete metric space $X$.

We can then prove the following
\bete\label{teo:attr1}
Let the assumptions of\/ {\rm Theorem~\ref{teo:smooth}} hold and suppose that
\begin{equation}
\label{diss}
\lambda_0 := \ess\sup_{x\in\Omega} \lambda (x) < 0.
\end{equation}
Then the dynamical system $(X,S(t))$ has a finite-dimensional
connected global attractor.
\ente
\noindent%
For readers' convenience, we report here below the argument of \cite{G}.
\begin{proof}
Consider $(\teta_{0i},\chi_{0i})\in X$, $i=1,2$,
set
\begin{equation*}
(\teta(t),\chi(t)) = ((\teta_1-\teta_2)(t),(\chi_1-\chi_2)(t))
\end{equation*}
where $(\teta_i(t),\chi_i(t))=S(t)(\teta_{0i},\chi_{0i})$ for $t\geq 0$, and
observe that
\beal{calorediff}
  & \teta_t+\chi_t+A\teta=0, \qquext{a.e.~in }\,(0,+\infty)\times\Omega,\\
  \label{phasediff}
  & \chi_t+J[\chi]+ f(\cdot,\chi_1) - f(\cdot,\chi_2)=\teta,
  \qquext{a.e.~in }\, (0,+\infty)\times\Omega.
\end{align}
Let us multiply equation \eqref{calorediff} by $ A\teta(t)$.
Integrating over $\Omega$, we get
\begin{equation*}
\frac1 2 \ddt\Vert\nabla\teta\Vert^2
    +\Vert A\teta\Vert^2=-(\chi_t, A\teta),
\end{equation*}
from which, using the Young and Poincar\'{e} inequalities, we derive
the estimate
\begin{equation*}
  \ddt\Vert\nabla\teta\Vert^2
    + \kappa_1\Vert \nabla\teta\Vert^2 \leq \Vert\chi_t\Vert^2.
\end{equation*}
and, by comparison in \eqref{phasediff}, we deduce
\begin{equation}
\label{base0}
\ddt\Vert\nabla\teta\Vert^2
    +\kappa_1\Vert \nabla\teta\Vert^2 \leq c_1\left (
    \Vert\chi\Vert^2 + \Vert J[\chi]\Vert^2
    +\Vert\teta\Vert^2\right),
\end{equation}
for all $t\geq 0$.
On the other hand, multiplying \eqref{phasediff} by $ \chi(t)$ and
integrating over $\Omega$ we find
\begin{equation*}
  \frac12\ddt \Vert\chi\Vert^2 + (f(\cdot,\chi_1) - f(\cdot,\chi_2),\chi)
    = - ( J[\chi],\chi) + (\teta, \chi).
\end{equation*}
Thus, on account of \eqref{defif}, \eqref{defif0} and \eqref{diss}, we get
\begin{equation*}
 \frac12\ddt \Vert\chi\Vert^2 + \lambda_0\Vert \chi\Vert^2
  \leq - ( J[\chi],\chi) + (\teta, \chi).
\end{equation*}
Then, Young's inequality gives
\begin{equation*}
\frac12\ddt \Vert\chi\Vert^2 + \frac{\lambda_0}{2}\Vert \chi\Vert^2
\leq c_{\lambda_0} \left(\Vert J[\chi]\Vert^2 + \Vert\teta\Vert^2\right).
\end{equation*}
Since $J$ is compact and self-adjoint on $H$,
we can find a finite-dimensional
projector $\Pi_{\lambda_0}$ such that
\begin{equation}
 \label{proj}
 \Vert J[v]\Vert^2 \leq \frac{\lambda_0}{4c_{\lambda_0}}\Vert v\Vert^2
  + c \Vert \Pi_{\lambda_0}[v]\Vert^2,
\end{equation}
for all $v\in H$. As a consequence we have
\begin{equation}
\label{base1}
\frac12\ddt \Vert\chi\Vert^2 + \frac{\lambda_0}{4}\Vert \chi\Vert^2
\leq c_{\lambda_0} \left(\Vert\Pi[\chi]\Vert^2 + \Vert\teta\Vert^2\right).
\end{equation}
Adding inequality \eqref{base0} multiplied by
$\mu=\frac{\lambda_0}{16c_1}$ to
\eqref{base1} and using the analogue of \eqref{proj}
to estimate the term $\| J[\chi] \|^2$ on the \rhs\
of \eqref{base0}, we infer
\begin{equation}\label{base2}
  \ddt \left(\mu\Vert\nabla\teta\Vert^2 + \frac12 \Vert\chi\Vert^2 \right)
   + \kappa_1\mu\Vert\nabla\teta\Vert^2 + \frac{\lambda_0}{8}
   \Vert \chi\Vert^2
  \leq c_{\lambda_0}\left(\Vert \Pi[\chi]\Vert^2
   + \Vert\teta\Vert^2\right).
\end{equation}
Therefore, from \eqref{base2}, we deduce
\begin{align} \label{base3}
\Vert\teta(t)\Vert_V^2 +\Vert \chi(t)\Vert^2 &\leq
 c_{\lambda_0}e^{-\kappa_{\lambda_0}t}
\left(\Vert\teta(0)\Vert_V^2 + \Vert\chi(0)\Vert^2\right)\\
\nonumber
&+ c_{\lambda_0}\int_{0}^t \,\left(\Vert\teta(\tau)\Vert^2
+ \Vert \Pi_{\lambda_0}[\chi(\tau)]\Vert^2 \right) \ditau,
\end{align}
for all $t\in [0,T]$ and any fixed $T>0$.

We now introduce the following pseudometric in $X$
$$
  \mathbf{d}_T ((\teta_{01},\chi_{01}),(\teta_{02},\chi_{02}) )
   = \left(\int_{0}^T\,\left(\Vert
   (\teta_1 -\teta_2)(\tau)\Vert^2 + \Vert \Pi_{\lambda_0}[(\chi_1-\chi_2)(\tau)]
   \Vert^2\right) \ditau\right)^{1/2}
$$
and we recall that a pseudometric is (pre)compact in $X$ (with
respect to the topology induced by the  $X$-metric) if any bounded
sequence in $X$ contains a Cauchy subsequence with respect to
$\mathbf{d}_T$ (see, for instance, \cite[Def.~2.57]{KM}).

Let us prove that $\mathbf{d}_T$ is precompact in $X$. Let
$\{(\teta_{0n},\chi_{0n})\}\subset X$ ($X$ is bounded) and set
$(\teta_n(t),\chi_n(t))=S(t)(\teta_{0n},\chi_{0n})$. Thanks to
\eqref{regoteta}, we have that $\{\teta_n\}$ is bounded in
$L^2(0,T;D(A))\cap H^1(0,T;H)$. Therefore it contains a
subsequence which strongly converges in $L^2(0,T,V)$. On the
other hand, we have that $\{\Pi_{\lambda_0}[\chi_n]\}$ is
bounded in $L^\infty(0,T;H)$. Also, by comparison in
\eqref{phase}, we deduce that $\{(\chi_n)_t\}$ is bounded in
$L^\infty(0,T;H)$. Therefore $\{(\Pi_{\lambda_0}[\chi_n])_t\}$
is bounded in $L^\infty(0,T;H)$ as well. Then
$\{\Pi_{\lambda_0}[\chi_{n}(\cdot)]\}$ contains a subsequence
which strongly converges in $L^2(0,T;H)$. Summing up
$\{(\teta_{0n},\chi_{0n})\}$ contains a Cauchy subsequence with
respect to $\mathbf{d}_T$.

From \eqref{base3}, we deduce that there exists $t^*>0$ such that
\begin{align*}
&\Vert S(t^*)(\teta_{01},\chi_{01}) - S(t^*)(\teta_{02},\chi_{02})
\Vert_X\\
&\leq
\frac{1}{2}\Vert(\teta_{01} - \teta_{02},\chi_{01} - \chi_{02})\Vert_X
+\bar{c}_{\lambda_0}
\mathbf{d}_{t^*} ((\teta_{01},\chi_{01}),
(\teta_{02},\chi_{02})).
\end{align*}
Hence $S(t)$ has a (connected) global attractor (see
\cite[Thm.~2.56, Prop.~2.59]{KM}) of finite fractal (i.e., box
counting) dimension (cfr. \cite[Thm.~2.8.1]{Ha}).
\end{proof}
\noindent%
In the case of singular potentials, we can first notice that,
by \eqref{regofchi} and the first \eqref{regochi2},
$t\mapsto \int_\Omega F_0(\chi(t))$ is absolutely continuous in $[0,T]$
due to \cite[Lemma~3.3, p.~73]{Br}. Hence,
taking $\calX$ defined by \eqref{defiXsing} as phase space, the
semigroup $S(t)$ defined as above takes $\calX$ to itself for all
$t\geq 0$. Using again \eqref{absorb} and on account of the
separation property \eqref{separ}, it is not difficult to realize that
there exists an absorbing set of the following form (note that
\eqref{conto23} still holds):
\begin{equation}\label{abssing}
  \mathcal{B}(R,\beta) := \big\{(\vartheta,\psi)\,:\,
   \Vert \vartheta\Vert_V\leq R,\; -1+\beta\leq \psi \leq 1-\beta~
   \text{a.e.~in~}\Omega\big\},
\end{equation}
for a suitable pair of constants $(R,\beta)\in
(0,+\infty)\times(0,1)$. Then, reasoning as above, we can suppose
that $\mathcal{B}(R,\beta)$ is invariant for $S(t)$ and we can
endow it with the $V\times H$-metric. The resulting complete metric
space $X$ is now our phase space and we have
\bete\label{teo:attr2}
Let the assumptions of\/ {\rm Theorem~\ref{teo:sing}} and \eqref{diss}
hold. Then, the dynamical system $(X,S(t))$ has a finite dimensional
connected global attractor.
\ente
\noindent
The proof goes as above.
\beos\label{asymptbehav}
 Assumption \eqref{diss} is crucial for
 investigating the asymptotic behavior (see, e.g., \cite[Thm.~1.2,
 (1.19) and Rem.~3]{FIRP}, \cite[Thm.~4.6 and (A4)]{BHZ} and
 \cite[Ass.~4, Sec.~2 and Lemma~3.1]{CCE2}). Also, we recall that if
 $f_0$ is real analytic then, following \cite{FIRP,GPS2},
 one can prove that the $\omega$-limit set of any
 pair $(\teta_0,\chi_0)$ in the energy space
 $\calX$ reduces to a singleton $\{(0,\chi_\infty)\}$, where
 $\chi_0$ solves the stationary problem
 $$
   J[\chi_\infty] + f_0(\chi_\infty) - \lambda(x) \chi_\infty= 0,
    \qquad\text{a.e.~in }\,\Omega.
 $$
\eddos

\bigskip\noindent
{\bf Acknowledgment.} This work was partially supported by the
Italian MIUR-PRIN Research Project 2008 ``Transizioni di fase,
isteresi e scale multiple''.




\end{document}